\documentclass[11pt]{article}
\usepackage{amsfonts}

\usepackage{CJK}
\usepackage{amsfonts,amssymb,amsmath,mathrsfs,multirow,booktabs}

\usepackage{color,latexsym,amsfonts}
 \newcommand{\qed}{\hfill\rule{2mm}{3mm}\vspace{4mm}}

 \newtheorem{theorem}{Theorem}[section]
 \newtheorem{lemma}[theorem]{Lemma}
 \newtheorem{corollary}[theorem]{Corollary}
 \newtheorem{proposition}[theorem]{Proposition}
 
 \newtheorem{Definition}[theorem]{Definition}
 \newtheorem{remark}[theorem]{Remark}
 \newtheorem{condition}[theorem]{Condition}

 \def\blemma{\begin{lemma}\sl{}\def\elemma{\end{lemma}}}
 \def\btheorem{\begin{theorem}\sl{}\def\etheorem{\end{theorem}}}

 \def\bcondition{\begin{condition}\sl{}\def\econdition{\end{condition}}}

 \def\beqlb{\begin{eqnarray}}\def\eeqlb{\end{eqnarray}}
 \def\beqnn{\begin{eqnarray*}}\def\eeqnn{\end{eqnarray*}}

 \def\mbb{\mathbb}\def\mbf{\mathbf}

 \def\va{\varepsilon}

 \def\<{\langle}\def\>{\rangle}

 \def\ar{&\!\!}

 \def\eqref#1{{\rm(\ref{#1})}}

 \def\proof{\noindent{\it Proof.~}}\def\qed{\hfill$\Box$\medskip}

\begin{document}
\noindent{(Draft: 2016/08/28)}
\bigskip

\centerline{\Large{\bf Maximum likelihood type estimation}}
\smallskip
\centerline{\Large{\bf for discretely observed CIR model}}
\smallskip
\centerline{\Large{\bf with small $\alpha$-stable
noises\footnote{Supported by NSFC (No.~11401012).
\textit{MSC (2010)}:
60G52; 62F12.
\textit{Keywords}: CIR model,
$\alpha$-stable, small noises, maximum likelihood type estimation.}}}
\bigskip

\centerline{Xu Yang\footnote{School of Mathematics and Information Science,
Beifang University of Nationalities, Yinchuan 750021, People's Republic of China.
Email: xuyang@mail.bnu.edu.cn}}

\bigskip

{\narrower{\narrower

\noindent{\bf Abstract.} A maximum likelihood type estimation of the
drift and volatility coefficient parameters in the CIR type model
driven by $\alpha$-stable noises is studied when the dispersion
parameter $\varepsilon\to0$ and the discrete observations frequency
$n\to\infty$ simultaneously.
\par}\par}

\section{Introduction}

\setcounter{equation}{0}

In mathematical finance, the classical Cox-Ingersoll-Ross (CIR)
model describes the evolution of interest rates. It specifies that
the instantaneous interest rate follows the stochastic differential
equation (SDE):
 \beqlb\label{1.0}
dx_\va(t)=(a_1'-a_2'x_\va(t))dt +a_3'\va\sqrt{x_\va(t)}dB(t),
 \eeqlb
where $\va,a_1',a_2',a_3'$ are strictly positive constants and
$\{B_t:t\ge0\}$ is a standard Brownian motion.

It is well-known that many financial processes exhibit discontinuous sample paths and
heavy tailed properties (e.g. certain moments are infinite).
These features cannot be captured by the CIR model. It is natural
to replace the driving Brownian motion by an $\alpha$-stable process;
see \cite{CW03} for the application of $\alpha$-stable processes in finance.
In this paper we are interested in the following stable driven CIR-type model:
 \beqlb\label{1.1}
dy_\va(t)=( a_1- a_2y_\va(t))dt + a_3\va y_\va(s-)^{1/q}dz_0(t),\qquad y_\va(0)=x_0\ge0,
 \eeqlb
where $a_1,a_3\ge0$, $q>0$, $
a_2\in\mbb{R}$ are constants, and $\{z_0(t):t\ge0\}$ is a spectrally
positive stable L\'evy process with index $\alpha\in(1,2)$ and
L\'evy measure $\mu(dz):=z^{-1-\alpha}1_{\{z>0\}}dz$.
By \cite[Corollary 4.3]{LMy10}, there is a pathwise
unique positive strong solution $\{y_\va(t):t\ge0\}$ to \eqref{1.1}
as $\frac1q+\frac1\alpha\ge1$. In the case of $q=\alpha$, the
solution is a particular form of the continuous-state branching
processes with immigration (see \cite[p.3]{Fu and Li}), which is
also called the stable CIR model (see \cite{LM13}).
If $a_1=a_2=0$, the solution can be treated as a critical branching process with
population dependent branching rate by \cite{Zhou16}.

Assume that the unknown quantity in \eqref{1.1} are the parameters
$a_1, a_2, a_3$. The type of data considered in this paper is
discrete observations at $n$ regularly spaced time points $t_k=k/n$
on the fixed interval $[0,1]$, that is $(y_\va(t_k))_{1\le k\le n}$.
The purpose of this paper is to study the maximum likelihood
estimator for the true value of $\textbf{a}:=( a_1, a_2, a_3)$ based
on these observations with small dispersion $\va$ and large sample
size $n$. To be precise, the type of asymptotics considered is when
$\va=\va_n$ goes to $0$ and $n$ goes to $\infty$ simultaneously. The
scheme of observations usually arises from some applied problems
such as the identification of a real deterministic dynamic system
with small random perturbations. We refer to \cite{UY04} for an application of small dispersion asymptotics to
contingent claim pricing.

The parameters estimation for discretely observed stochastic
processes driven by small Brownian motion has been studied by
several authors; see e.g. \cite{SU03,GS09,Uch04,Uch08}. The
asymptotics distributions of the estimators based on a Gaussian approximation to the transition density (see  \cite{Kessler97}) are normal under certain
conditions on $\va=\va_n$ and $n$; see e.g. \cite{SU03}, where $n\to\infty$ and $\lim_{n\to\infty}(\va
\sqrt{n})^{-1}<\infty$.

Recently, a number of papers have been devoted to small volatility
asymptotics for the parameter estimation in the models driven by small
L\'evy noises. When the coefficient of the L\'evy jump term is
constant, drift parameter estimation of discretely observed
L\'evy driven SDEs has been studied by many authors; see e.g. \cite{Long 09,M10,Long13}. For the SDE \eqref{1.1}, where the jumps are state-dependent
and the jump term is \textit{non-Lipschitz}, the asymptotics
properties of the conditional least squares estimators and the
weighted conditional least squares estimators of the drift
parameters ($a_1,a_2$) were given in \cite{LM13} based on low frequency observations, and the asymptotics behavior of the
least squares estimator of the parameter $a_1$ (or $a_2$) was
established in \cite{MaYang2014} under high frequency observations and small
dispersion.

In this paper we employ a maximum likelihood type method to obtain
an estimator for the parameter $\textbf{a}=( a_1, a_2, a_3)$ in
\eqref{1.1}. To overcome the difficulty that the joint density of the sample
$\{y_\va(t):t\in[0,1]\}$ is not tractable, we deal with
by \textit{using stable distributions to approximate the density}. It follows
from \eqref{1.1} that
$$y_\va(t_k)=y_\va(t_{k-1})+ a_1\Delta t_k -
a_2\int_{t_{k-1}}^{t_k}y_\va(s)ds +  a_3\va
\int_{t_{k-1}}^{t_k}y_\va(s-)^{1/q}dz_0(s),$$
where $\Delta t_k=t_k-t_{k-1}=1/n$. Then one can use the Euler
scheme (see e.g. \cite{Jacod04}, which studied for a SDE driven
by a L\'evy process) to get the approximation
$$y_\va(t_k)\approx y_\va(t_{k-1})+ a_1\Delta t_k -
a_2y_\va(t_{k-1})\Delta t_k + a_3\va \Delta t_k^{1/\alpha}
y_\va(t_{k-1})^{1/q}z_k,$$
where $z_1,z_2,\cdots,z_n$ are independent stable random variables
with the same distribution as $z_0(1)$.
So as $\Delta t_k$, the distance between
observations, is small, it may suggest that,
conditioned on $y_\va(t_{k-1})$, the distribution of this random variable
$y_\va(t_k)-y_\va(t_{k-1})-a_1\Delta t_k + a_2y_\va(t_{k-1})\Delta
t_k$ may be close to that of
$a_3\va \Delta t_k^{1/\alpha} y_\va(t_{k-1})^{1/q}z_k$ in ceratin sense.
Inspired by this, we can define a \textit{likelihood type function} of
$(y_\va(t_k))_{0\le k\le n}$ by
 \beqnn
L_{\va,n}(\textbf{a}) :=\prod_{k=1}^n [a_3\va
n^{-1/\alpha}y_\va(t_{k-1})^{1/q}]^{-1}p(Y_{\va,n,k}(\textbf{a})),
 \eeqnn
where $p(x)$ is the density function of $z_0(1)$ and
 \beqlb\label{1.2}
Y_{\va,n,k}( \textbf{a})
 \ar:=\ar
[y_\va(t_k)-y_\va(t_{k-1})- a_1\Delta t_k + a_2
y_\va(t_{k-1})\Delta t_k]\cdot[a_3\va \Delta t_k^{1/\alpha}
y_\va(t_{k-1})^{1/q}]^{-1} \cr
 \ar=\ar
\big[y_\va(t_k)-y_\va(t_{k-1})- a_1/n + a_2
y_\va(t_{k-1})/n\big]\cdot\big[a_3\va
n^{-1/\alpha}y_\va(t_{k-1})^{1/q}\big]^{-1}.
 \eeqlb
This likelihood type function may be a bit like the joint density of $(y_\va(t_k))_{0\le k\le n}$ as $\Delta t_k$ is small enough.
Now we define the \textit{log likelihood type function} of $(y_\va(t_k))_{0\le k\le n}$ by
 \beqnn
\tilde{U}_{\va,n}(\textbf{a})
 :=
\log L_{\va,n}(\textbf{a})
 =
\sum_{k=1}^n\log p(Y_{\va,n,k}(\textbf{a})) -n\log  a_3
-q^{-1}\sum_{k=1}^n\log y_\va(t_{k-1})-n\log(\va n^{-1/\alpha}).
 \eeqnn
Let
$\hat{\textbf{a}}_{\va,n}:=(\hat{ a}_{1,\va,n},\hat{
a}_{2,\va,n},\hat{ a}_{3,\va,n})$ be the \textit{maximum likelihood type estimator}
defined by $
\tilde{U}_{\va,n}(\hat{\textbf{a}}_{\va,n})= \sup_{
\textbf{a}\in{\bf\bar{A}}} \tilde{U}_{\va,n}(\textbf{a})$,
where ${\bf\bar{A}}$ is the closure of an open set defined in
Section 2. Such approximation is usually called an Euler-Maruyama approximation in the classical CIR model defined by \eqref{1.0}; see e.g. \cite[Section
9.1]{KP92} and \cite[Section 4.2.2]{R04}.

It is obvious that
$\hat{\textbf{a}}_{\va,n}$ is also a maximum point of $U_{\va,n}$
defined by
 \beqlb\label{1.a}
U_{\va,n}(\textbf{a}):=\sum_{k=1}^n\log p(Y_{\va,n,k}(\textbf{a}))
-n\log  a_3,\qquad U_{\va,n}(\hat{\textbf{a}}_{\va,n})= \sup_{
\textbf{a}\in{\bf \bar{A}}} U_{\va,n}(\textbf{a}).
 \eeqlb
Our \textit{main result} of this paper, Theorem \ref{t1.1}, gives a
consistent, asymptotically normal and asymptotically efficient
estimator $\hat{\textbf{a}}_{\va,n}$ of $\textbf{a}$ under the
conditions $\va=\va_n\to0$, $n\to\infty$ and $\lim_{n\to\infty}(\va
n^{1/\alpha-1})^{-1}<\infty$, which is consistent with the
corresponding assertion in \cite[Theorem 1]{SU03} if $\alpha=q=2$.
The proof is established in Section 3. An auxiliary lemma
and the proof of Lemma \ref{t1.2} are presented in Section 4.

\section{Main result}

\setcounter{equation}{0}

Before stating the main result of this paper, we give some
notations. We always assume that all random elements are defined on
a filtered complete probability space $(\Omega, \mathscr{F}, (
\mathscr{F}_t)_{t\in[0,1]},\mathbf{P})$ satisfying the usual
hypotheses. Let $C(\mbb{R})$ be the space of continuous functions on
$\mbb{R}$. For $\mbb{R}_+^2:=[0,\infty)\times[0,\infty)$ define
$C(\mbb{R}_+^2)$ similarly. For $f,g\in C(\mbb{R})$ write $\<f,g\>=\int_{\mbb{R}}f(x)g(x)dx$. For any integer $n\ge 1$ let
$C^n(\mbb{R})$ be the subset of $C(\mbb{R})$ with continuous
derivatives up to the $n$th order. Set $\|x\| =
\sup_{t\in[0,1]}|x(t)|$. We use $`` \overset{p}{\longrightarrow}"$
and $`` \overset{d}{\longrightarrow}"$ to denote the convergence of
random variables in probability and in distribution, respectively.
Let $\textbf{a}$ be the parameter and $\bar{{\bf
a}}:=(\bar{a}_1,\bar{a}_2,\bar{a}_3)$ the true value of
$\textbf{a}$. For $t\in[0,1]$ define $y_0(t)=x_0e^{-\bar{ a}_2t}+\bar{
a}_1\int_0^te^{-\bar{ a}_2(t-s)}ds$. Put
 \beqlb\label{1.b}
U(\textbf{a})=\int_0^1dt \int_{\mbb{R}}p(x)\log p(Y_0(
\textbf{a},t,x))dx -\log a_3,
 \eeqlb
where $
Y_0( \textbf{a},t,x):=m_0Y(\textbf{a},t)+\bar{ a}_3 a_3^{-1}x$ and
$Y(\textbf{a},t):=(\bar{ a}_1- a_1)a_3^{-1}y_0(t)^{-1/q} +(
a_2-\bar{ a}_2)a_3^{-1}y_0(t)^{1-1/q}$.
Define the matrices ${\bf V}({\bf a}):=(\frac{\partial^2 U({\bf
a})}{\partial a_i\partial  a_j})$ and
 \beqnn
{\bf\Sigma}:=\left(\begin{matrix}
 v_1m^{0,2} \qquad
 -v_1m^{1,2} \qquad
 v_2m^{0,1} \\
 -v_1m^{1,2} \qquad
 v_1m^{2,2} \qquad
 -v_2m^{1,1} \\
 v_2m^{0,1} \qquad
 -v_2m^{1,1} \qquad\quad
 v_3~~
\end{matrix}\right),
 \eeqnn
where
$v_1:=\int_{\mbb{R}}|p'(x)|^2/p(x)dx$,
$v_2:=\int_{\mbb{R}}x|p'(x)|^2/p(x)dx$,
$v_3:=\int_{\mbb{R}}x^2|p'(x)|^2/p(x)dx-1$, and
$m^{i,j}:=\int_0^1y_0(t)^{i-\frac{j}{q}}dt$ for $i,j\ge0$.
We give the conditions on the initial value $x_0=y_\va(0)$ and the true value of the parameters.
 \bcondition\label{c1.1} Neither of the following conditions
hold: (i) $x_0=\bar{a}_1/\bar{a}_2$ and $\bar{a}_2\neq0$; (ii)
$\bar{a}_1=0$ and $\bar{a}_2=0$; (iii) $x_0=\bar{a}_1=0$.
 \econdition
Observe that $y_0(t)>0$ for all $t>0$ under Condition \ref{c1.1}.
Since $\lim_{\va\to0}\|y_\va-y_0\|=0~\mbf{P}$-a.s. by \cite[Proposition 3.2]{MaYang2014}, $y_\va(t)>0$ for all $t\in[0,1]$ as $\va$ small enough.
This makes sure that \eqref{1.2} is well defined.
It is easy to see that $\frac{\partial U(\bar{{\bf a}})}{\partial
a_i}=0$ ($i=1,2,3$) and ${\bf V}({\bf
\bar{a}})=-\bar{a}_3^{-2}{\bf\Sigma}$. Then ${\bf \bar{a}}$ is a
local maximum point of $U(\textbf{a})$ under Condition \ref{c1.1} by
Lemma \ref{t1.6} in Appendix. In the following we state the conditions
on the domain ${\bf A}$ and the relationship between $n$ and $\va$.
 \bcondition\label{c1.2}
(i) Let ${\bf A}$ be an open bounded convex subset of $[0,\infty)\times
\mbb{R}\times[0,\infty)$ and ${\bf\bar{ A}}$ denote its closure
set. Suppose that ${\bf\bar{
A}}\cap(\mbb{R}^2\times\{0\})=\emptyset$ and ${\bf \bar{a}}\in {\bf
A}$ is the only maximum point of $U(\textbf{a})$ on ${\bf\bar{ A}}$.
(ii) Suppose that $\va:=\va_n$ and
$\lim_{n\to\infty}m_{\va,n}=m_0<\infty$, where
$m_{\va,n}:=\va^{-1}n^{\frac1\alpha-1}$.
 \econdition
 \btheorem\label{t1.1} Assume that Conditions \ref{c1.1}--\ref{c1.2}
hold. Then as $n\to\infty$,
\begin{gather}
\hat{ {\bf a}}_{\va,n}\overset{p}{\longrightarrow} \bar{ {\bf a}},\label{1.4}\\
\textbf{S}_{\va,n}:=\big( v_{\va,n}(\hat{ a}_{1,\va,n}-\bar{ a}_1),
v_{\va,n}(\hat{ a}_{2,\va,n}-\bar{ a}_2), \sqrt{n}(\hat{
a}_{3,\va,n}-\bar{ a}_3) \big)\overset{d}{\longrightarrow}
N(\textbf{0},\bar{a}_3^2{\bf\Sigma}^{-1}),\label{1.5}
\end{gather}
where
$v_{\va,n}:=m_{\va,n}\sqrt{n}=\va^{-1}n^{\frac1\alpha-\frac12}$ and
$\textbf{0}:=(0,0,0)$.
 \etheorem



\section{Proof of Theorem \ref{t1.1}}

\setcounter{equation}{0}

It
follows from \eqref{1.1} that
 \beqlb\label{2.5}
y_\va(t_k)-y_\va(t_{k-1}) =\frac{\bar{ a}_1}{n}-\bar{
a}_2\int_{t_{k-1}}^{t_k}y_\va(s)ds +\va\bar{
a}_3\int_{t_{k-1}}^{t_k} y_\va(s-)^{1/q}dz_0(s).
 \eeqlb
Together with \eqref{1.2} one derives that for $ \textbf{a}=( a_1,
a_2, a_3)$,
$n,k\ge1$ and $\va>0$,
 \beqlb\label{1.3}
Y_{\va,n,k}(\textbf{a})
 \ar=\ar
(\bar{ a}_1- a_1)
a_3^{-1}m_{\va,n}n\int_{t_{k-1}}^{t_k}y_\va([ns]/n)^{-1/q}ds\cr
 \ar\ar
+ a_3^{-1}m_{\va,n} n\int_{t_{k-1}}^{t_k}\big[
a_2y_\va([ns]/n)^{1-1/q} -\bar{
a}_2y_\va(s)y_\va([ns]/n)^{-1/q}\big]ds\cr
 \ar\ar
+\bar{a}_3 a_3^{-1}n^{1/\alpha}\int_{t_{k-1}}^{t_k}
\big[y_\va(s-)y_\va([ns]/n)^{-1}\big]^{1/q}dz_0(s) \cr
 \ar=:\ar
(\bar{ a}_1- a_1) a_3^{-1}m_{\va,n}N_{\va,n,k} +
a_3^{-1}m_{\va,n}\bar{M}_{\va,n,k}( a_2) +\bar{ a}_3
a_3^{-1}K_{\va,n,k},
 \eeqlb
where $[x]$ denotes the largest integer not greater than $x$, and
$M_{\va,n,k}=n\int_{t_{k-1}}^{t_k}y_\va([ns]/n)^{1-1/q}ds$
and $K_{n,k} = n^{1/\alpha}[z(t_k)-z(t_{k-1})]$.
Before showing the proof of Theorem \ref{t1.1}, we state the
following lemma, which will be proved in Appendix.
 \blemma\label{t1.2} Suppose that Conditions \ref{c1.1}--\ref{c1.2} hold.
 Let $H\in C^1(\mbb{R})$ and
$B\in C(\mbb{R}_+^2)$ satisfy
 \beqlb\label{2.30}
\sup_{x\ge0}|H'(x)|+\sup_{x<0}|(-x)^{-\gamma}H'(x)|<\infty
 \eeqlb
and
$|B(x_1,y_1)-B(x_2,y_2)|\le C_k[|x_1-x_2|+|y_1-y_2|]$ for all
$x_1,x_2,y_1,y_2\in[0,k]$ and $k\ge1$, where $\gamma>1$ and $C_k>0$ are constants. Then
as $n\to\infty$,
 \beqlb\label{2.34}
\sup_{ \textbf{a}\in{\bf \bar{A}}}\Big|\frac1n\sum_{k=1}^n
H(Y_{\va,n,k}(\textbf{a}))B(M_{\va,n,k},N_{\va,n,k}) -\int_0^1B_tdt
\int_{\mbb{R}}p(x)H(Y_0( \textbf{a},t,x))dx\Big|
\overset{p}{\longrightarrow}0.
 \eeqlb
Moreover, if $\bar{H}\in C^1(\mbb{R})$ satisfying \eqref{2.30} with
$H$ replaced by $\bar{H}$ and $\<\bar{H},p\>=\<H,p\>=0$, then
 \beqlb\label{2.15}
n^{-1/2}\sum_{k=1}^n
\Big[H(Y_{\va,n,k}(\bar{\textbf{a}}))B(M_{\va,n,k},N_{\va,n,k})
+\bar{H}(Y_{\va,n,k}(\bar{\textbf{a}}))\Big]
\overset{d}{\longrightarrow}N(0,\eta_0^2)
 \eeqlb
as $n\to\infty$, where
$\eta_0^2:=\int_0^1
\<(B_t H+\bar{H})^2,p\>dt$ and
$B_t:=B(y_0(t)^{1-1/q},y_0(t)^{-1/q})$.
 \elemma

 \blemma\label{t1.3} For each $k\ge0$, there are constants
$c_k,c_k'>0$ so that
$p^{(k)}(x)\sim c_kx^{-\alpha-1-k}$ and $p^{(k)}(-x)\sim
c_k'\xi^{\frac{2-\alpha+2k}{2\alpha}}e^{-\xi}$
as $x\to\infty$, where
$\xi=(\alpha-1)(x/\alpha)^{\alpha/(\alpha-1)}$.
 \elemma
 \proof The proofs of these two assertions follow immediately from
the arguments in \cite[Theorem 2.5.1]{Zolotarev86} and \cite[Theorem
2.5.2]{Zolotarev86}, respectively. \qed

For $x\in\mbb{R}$ define
$H_0(x)
=
p'(x)/p(x)$,
$H_1(x)= \big[p''(x)p(x)-|p'(x)|^2\big]/p(x)^2$,
$H_2(x) = xH_1(x)+H_0(x)$ and $H_3(x) =x^2H_1(x)+2xH_0(x)+1$.
For $i,j=1,2,3$ let $U_{\va,n}^i(\textbf{a})=\frac{\partial
U_{\va,n}(\textbf{a})}{\partial
 a_i}$ and $U_{\va,n}^{i,j}(\textbf{a})=\frac{\partial^2
U_{\va,n}(\textbf{a})}{\partial  a_i\partial  a_j}$. For $1\le
i_1,j_1\le2$ and $(i_2,j_2)\in\{(1,3),(2,3),(3,1),(3,2)\}$ let
$V_{\va,n}^{i_1,j_1}(\textbf{a})=
v_{\va,n}^{-2}U_{\va,n}^{i_1,j_1}(\textbf{a})$,
$V_{\va,n}^{i_2,j_2}(\textbf{a})=
v_{\va,n}^{-1}n^{-\frac12}U_{\va,n}^{i_2,j_2}(\textbf{a})$
and
 \beqnn
 \ar\ar
V^{i_1,j_1}(\textbf{a})=
(-1)^{i_1+j_1}a_3^{-2}\int_0^1y_0(t)^{i_1+j_1-2-\frac2q}
\<H_1(Y_0(\textbf{a},t,\cdot),p\>dt
, \cr
 \ar\ar
V^{i_2,j_2}(\textbf{a})=
(-1)^{i_2+j_2}a_3^{-2}\int_0^1y_0(t)^{i_2+j_2-4-\frac1q}
\<H_2(Y_0(\textbf{a},t,\cdot),p\>dt.
 \eeqnn
Put $
V_{\va,n}^{3,3}(\textbf{a})= n^{-1}U_{\va,n}^{3,3}(\textbf{a})$ and
$V^{3,3}(\textbf{a}):=
a_3^{-2}\<H_3(Y_0(\textbf{a},t,\cdot),p\>$.
Then ${\bf V}({\bf a})=(V^{i,j}(\textbf{a}))$. Define the matrix
$\textbf{V}_{\va,n}(\textbf{a})=(V_{\va,n}^{i,j}(\textbf{a}))$. Set
${\bf\Lambda}_{\va,n}=(
v_{\va,n}^{-1}U_{\va,n}^1(\bar{\textbf{a}}),
v_{\va,n}^{-1}U_{\va,n}^2(\bar{\textbf{a}}),
n^{-\frac12}U_{\va,n}^3(\bar{\textbf{a}}))$.

 \blemma\label{t1.5}
Suppose that Conditions \ref{c1.1}--\ref{c1.2} hold. Then as
$n\to\infty$,
\beqnn
{\bf\Lambda}_{\va,n}\overset{d}{\longrightarrow}
N(\textbf{0},\bar{a}_3^{-2}{\bf\Sigma}) \mbox{\quad and \quad}
\sup_{ \textbf{a}\in {\bf
\bar{A}}}|\textbf{V}_{\va,n}(\textbf{a})-\textbf{V}(\textbf{a})|
\overset{p}{\longrightarrow}0.
\eeqnn
 \elemma
 \proof
It is easy to see that
 \beqnn
v_{\va,n}^{-1}U_{\va,n}^1(\bar{\textbf{a}})
 \ar=\ar
-\bar{
a}_3^{-1}n^{-\frac12}
\sum_{k=1}^nH_0(Y_{\va,n,k}(\bar{\textbf{a}}))N_{\va,n,k}, \cr
v_{\va,n}^{-1}U_{\va,n}^2(\bar{\textbf{a}})
 \ar=\ar
\bar{
a}_3^{-1}n^{-\frac12}
\sum_{k=1}^nH_0(Y_{\va,n,k}(\bar{\textbf{a}}))M_{\va,n,k},\cr
n^{-\frac12}U_{\va,n}^3(\bar{\textbf{a}})
 \ar=\ar
-\bar{
a}_3^{-1}n^{-\frac12}
\sum_{k=1}^n\big[H_0(Y_{\va,n,k}(\bar{\textbf{a}}))Y_{\va,n,k}(\bar{\textbf{a}})+1\big].
 \eeqnn
Observe that $\<H_0,p\>=0$ and
$\int_{\mbb{R}}[xH_0(x)+1]p(x)dx=0$. Then by Lemmas \ref{t1.2} and
\ref{t1.3},
 \beqnn
x_1v_{\va,n}^{-1}U_{\va,n}^1(\bar{\textbf{a}})+x_2v_{\va,n}^{-1}U_{\va,n}^2(\bar{\textbf{a}})
+x_3n^{-\frac12}U_{\va,n}^3(\bar{\textbf{a}})
\overset{d}{\longrightarrow}N(0,\eta(x_1,x_2,x_3)^2)
 \eeqnn
for all $x_i\in\mbb{R}$ ($i=1,2,3$) as $n\to\infty$, where
 \beqnn
\eta(x_1,x_2,x_3)^2 :=
\bar{ a}_3^{-2}\int_0^1dt
\int_{\mbb{R}}\big[(x_1y_0(t)^{-\frac1q}
-x_2y_0(t)^{1-\frac1q})p'(x)+x_3(xp'(x)+p(x))\big]^2/p(x)dx.
 \eeqnn
Then the first assertion follows from the Cram\'er-Wold theorem.

It is elementary to see that for $1\le i_1,j_1\le2$ and
$(i_2,j_2)\in\{(1,3),(2,3),(3,1),(3,2)\}$,
 \beqnn
V_{\va,n}^{i_1,j_1}(\textbf{a})
 \ar=\ar
(-1)^{i_1+j_1}a_3^{-2}n^{-1}\sum_{k=1}^nH_1(Y_{\va,n,k}(\textbf{a}))
M_{\va,n,k}^{i_1+j_1-2}N_{\va,n,k}^{4-i_1-j_1}, \cr
V_{\va,n}^{i_2,j_2}(\textbf{a})
 \ar=\ar
(-1)^{i_2+j_2}a_3^{-2}n^{-1}\sum_{k=1}^n
H_2(Y_{\va,n,k}(\textbf{a}))
M_{\va,n,k}^{i_2+j_2-4}N_{\va,n,k}^{5-i_2-j_2}
 \eeqnn
and $V_{\va,n}^{3,3}(\textbf{a})
=\frac{a_3^{-2}}{n}\sum_{k=1}^nH_3(Y_{\va,n,k}(\textbf{a}))$. It
follows from Lemma \ref{t1.3} that $H_i\in C^1(\mbb{R})$ satisfies
\eqref{2.30} with $H$ replaced by $H_i$ for $i=1,2,3$. Therefore, by
Lemma \ref{t1.2} we know that
$\sup_{\textbf{a}\in {\bf \bar{A}}}|V_{\va,n}^{i,j}(\textbf{a})
-V^{i,j}(\textbf{a})|\overset{p}{\longrightarrow}0$
for each $1\le i,j\le3$, which derives the last assertion.
 \qed

\noindent{\it Proof of Theorem \ref{t1.1}.} The proof is a
modification of that of \cite[Theorem 1]{SU03}. We give some details
in the following. Suppose that there exists a subsequence
$(\va_{n_k},n_k)$ so that ${\bf\hat{ a}}_{\va_{n_k},n_k}$ tends to a
limit ${\bf\check{ a}}=(\check{ a}_1,\check{ a}_2,\check{ a}_3)$.
Taking $H(x)=\log p(x)$ in Lemma \ref{t1.2}, we get
 \beqlb\label{2.29}
\sup_{ \textbf{a}\in{\bf \bar{A}}}|n^{-1} U_{\va,n}(\textbf{a})
-U(\textbf{a})| \overset{p}{\longrightarrow}0,
 \eeqlb
where $U(\textbf{a})$ is defined in \eqref{1.b}.
By \eqref{1.a}, for each $k\ge1$ we get
$\frac1n_k U_{\va_{n_k},n_k}({\bf\bar{ a}})\le \frac1n_k
U_{\va_{n_k},n_k}(\hat{{\bf a}}_{\va_{n_k},n_k})$.
Letting $k\to\infty$, by \eqref{2.29}, we have $U({\bf\bar{ a}}) \le
U(\check{\textbf{a}})$. On the other hand, ${\bf\bar{ a}}$ is the
only maximum point of $U(\textbf{a})$ by Condition \ref{c1.2}(i), thus
$\check{\textbf{a}}={\bf\bar{ a}}$. This proves \eqref{1.4}.

By Taylor's formula, $\textbf{S}_{\va,n}\textbf{D}_{\va,n}={\bf\Lambda}_{\va,n},$
where $\textbf{D}_{\va,n} =\int_0^1\textbf{V}_{\va,n}({\bf\bar{
a}}+u(\hat{ {\bf a}}_{\va,n}-{\bf\bar{ a}}))du$. Then by Lemma
\ref{t1.5} and the fact ${\bf V}({\bf
\bar{a}})=-\bar{a}_3^{-2}{\bf\Sigma}$, one obtains \eqref{1.5} by using
the same argument in the corresponding proof of \cite[Theorem
1]{SU03}.
 \qed

\section{Appendix}
\setcounter{equation}{0}

 \blemma\label{t1.6}
Suppose that Condition \ref{c1.1} holds. Then ${\bf\Sigma}$ is a
positive definite matrix.
 \elemma
 \proof
By the H\"older inequality, the determinant
 \beqlb\label{2.33}
\left|\begin{matrix}
 v_1m^{0,2} \quad
 -v_1m^{1,2}  \\
 -v_1m^{1,2} \quad
 v_1m^{2,2}
\end{matrix}\right|
=v_1^2\big[m^{0,2}m^{2,2}-(m^{1,2})^2\big]>0
 \eeqlb
and
$1=|\int_{\mbb{R}}\frac{xp'(x)}{p(x)^{1/2}}p(x)^{\frac12}dx|^2
<\int_{\mbb{R}}[xp'(x)]^2/p(x)dx$,
which implies $v_3>0$. It is obvious that
 \beqlb\label{2.300}
|{\bf\Sigma}|v_3=\big[v_1 v_3m^{2,2}-(v_2m^{1,1})^2\big]\big[v_1
v_3m^{0,2}-(v_2m^{0,1})^2\big]-\big[v_1
v_3m^{1,2}-v_2^2m^{0,1}m^{1,1}\big]^2.
 \eeqlb
Since $\<1,p'\>=0$ and $\int_{\mbb{R}}xp'(x)dx=-1$, by
the H\"older inequality again we get
 \beqnn
v_2^2
 \ar=\ar
\Big|\int_{\mbb{R}}\frac{p'(x)}{p(x)}[xp'(x)+p(x)]dx\Big|^2 <
\int_{\mbb{R}}\frac{|p'(x)|^2}{p(x)}dx\int_{\mbb{R}}\frac{[xp'(x)+p(x)]^2}{p(x)}dx
\cr
 \ar=\ar
\int_{\mbb{R}}\frac{|p'(x)|^2}{p(x)}dx\Big[\int_{\mbb{R}}\frac{[xp'(x)]^2}{p(x)}dx
+2\int_{\mbb{R}}xp'(x)dx+1\Big]=v_1v_3.
 \eeqnn
Note that
$|\int_0^1[y_0(t)^{1-1/q}z+y_0(t)^{-1/q}]dt|^2 <
\int_0^1[y_0(t)^{1-1/q}z+y_0(t)^{-1/q}]^2dt$ for each $z\in\mbb{R}$.
Then
 \beqnn
 \ar\ar
\big[v_1 v_3m^{2,2}-(v_2m^{1,1})^2\big]z^2+2\big[v_1
v_3m^{1,2}-v_2^2m^{0,1}m^{1,1}\big]z+\big[v_1
v_3m^{0,2}-(v_2m^{0,1})^2\big] \cr
 \ar\ar
=
v_1v_3\big[m^{2,2}z^2+2zm^{1,2}+m^{0,2}\big]-v_2^2\big[(m^{1,1}z)^2+2zm^{1,1}m^{0,1}+(m^{0,1})^2\big]
\cr \ar\ar
 =
v_1v_3\int_0^1\big[y_0(t)^{1-1/q}z+y_0(t)^{-1/q}\big]^2dt
-v_2^2\Big|\int_0^1\big[y_0(t)^{1-1/q}z+y_0(t)^{-1/q}\big]dt\Big|^2>0
 \eeqnn
for all $z\in\mbb{R}$. It follows from \eqref{2.300} that $|{\bf\Sigma}|v_3>0$,
which implies $|{\bf\Sigma}|>0$. Together with \eqref{2.33} one gets
the desired result. \qed

\noindent{\it Proof of Lemma \ref{t1.2}.} In the following $C$ is a constant whose value might change from place to
place and does not depend on $\va,n,k,t$ and $\textbf{a}$. For
$\textbf{a}\in{\bf \bar{A}}$, $n,k\ge1$, $\va>0$ and $t\in[0,1]$ we
put
$B_{\va,n,k}:=B(M_{\va,n,k},N_{\va,n,k})$, 
$Y_{n,k}(\textbf{a},t)
 :=
m_0 Y(\textbf{a},t)+\bar{ a}_3 a_3^{-1}K_{n,k}$ and
$\bar{Y}_{n,k}(\textbf{a},t,\va):= m_{\va,n}Y(\textbf{a},t) +\bar{
a}_3 a_3^{-1}K_{n,k}$.
For $0<\zeta<\inf_{t\in[0,1]} y_0(t)$ and $\varsigma>\|y_0\|$ define
$A_{\va,\zeta}=\{\inf_{t\in[0,1]}y_\va(t)\ge\zeta\}$ and
$B_{\va,\varsigma}=\{\sup_{t\in[0,1]}y_\va(t) \le\varsigma\}$.
It follows from \cite[Lemma 3.5]{MaYang2014} that
 \beqlb\label{2.7}
\mbf{P}\{A_{\va,\zeta}^c\} +\mbf{P}\{B_{\va,\varsigma}^c\}\le
C\va^\alpha.
 \eeqlb
Define $U_{\va,\varsigma,\zeta}=A_{\va,\zeta}\cap
B_{\va,\varsigma}$. We divide the rest of proof into seven steps.

{\bf Step 1.} First we show: For
each $\gamma'>1$ and large enough $n_0>0$,
 \begin{gather}
\sup_{n,k\ge1,~\va>0} \mbf{E}\Big\{\sup_{\textbf{a}\in {\bf
\bar{A}}}[-Y_{\va,n,k}(\textbf{a})]^{\gamma'}1_{\{Y_{\va,n,k}(\textbf{a})<-n_0,
U_{\va,\varsigma,\zeta}\}}\Big\}<\infty,\label{2.2}\\
\sup_{n,k\ge1} \mbf{E}\big\{[-K_{n,k}]^{\gamma'}
1_{\{K_{n,k}<0\}}\big\}<\infty.\label{2.1}
 \end{gather}
By \cite[Theorem 1.4]{KSh02}, for each $\va>0$ and $n,k\ge1$ there
is a stable process $\{z_{\va,n,k}(t):t\ge0\}$ with the
same finite dimension distribution as $\{z_0(t):t\ge0\}$ so that
$K_{\va,n,k}=z_{\va,n,k}(T_{\va,n,k})$, where
$T_{\va,n,k}=n\int_{t_{k-1}}^{t_k}
[y_\va(s)y_\va([ns]/n)^{-1}]^{\alpha/q}ds$.
Observe that $T_{\va,n,k}\le (\varsigma\zeta^{-1})^{\alpha/q}$ on
$A_{\va,\zeta}\cap B_{\va,\varsigma}$. It follows from \cite[Lemma
2.4]{FMW10} that for each $x>0$,
 \beqlb\label{2.22}
\mbf{P}\big\{K_{\va,n,k}\le -x,
U_{\va,\varsigma,\zeta}\big\} \le \mbf{P}\big\{\inf_{t\le
\varsigma\zeta^{-1}} z_{\va,n,k}(t)\le -x\big\}
\le\exp\big\{-\tilde{c}_0x^{\alpha/(\alpha-1)}\big\},
 \eeqlb
where $\tilde{c}_0:=[(\alpha-1)/\alpha]^{\alpha/(\alpha-1)}
[\zeta\varsigma^{-1}]^{\frac{\alpha}{q(\alpha-1)}}$. One can also
see that
$\bar{M}_{\va,n,k}( a_2)\le | a_2|\varsigma^{1-1/q}+|\bar{
a}_2|\varsigma\zeta^{-1/q}$ and $N_{\va,n,k}\le \zeta^{-1/q}$
on $U_{\va,\varsigma,\zeta}$, which implies
$|Y_{\va,n,k}({\bf a})|^{\gamma'}\le \tilde{c}_1 +|2\bar{ a}_3
a_3^{-1}K_{\va,n,k}|^{\gamma'}$ on $U_{\va,\varsigma,\zeta}$ with
$$\tilde{c}_1:=\sup_{ {\bf a}\in {\bf \bar{A}},\,\va>0,\,n\ge1}
|2m_{\va,n} a_3^{-1}[| a_2|\varsigma^{1-1/q} +|\bar{
a}_2|\varsigma\zeta^{-1/q}+ | a_1-\bar{
a}_1|\zeta^{-1/q}]|^{\gamma'}.$$
Thus
 \beqlb\label{2.8}
 \ar\ar
\mbf{E}\Big\{\sup_{{\bf a}\in {\bf \bar{A}}}[-Y_{\va,n,k}({\bf
a})]^{\gamma'}1_{\{Y_{\va,n,k}({\bf a})<-n_0,
U_{\va,\varsigma,\zeta}\}}\Big\}
\le
C
\mbf{E}\Big\{[-K_{\va,n,k}]^{\gamma'}1_{\{K_{\va,n,k}<0,
U_{\va,\varsigma,\zeta}\}}\Big\} + \tilde{c}_1\cr
 \ar\ar\qquad=
C\int_0^\infty
t^{\gamma'-1}\mbf{P}\{K_{\va,n,k}<-t,
U_{\va,\varsigma,\zeta}\}dt + \tilde{c}_1
 \eeqlb
for large enough $n_0$. Together with \eqref{2.22} implies \eqref{2.2}. Similarly, one can
also get \eqref{2.1}.

{\bf Step 2.} In this step we show that for $\delta\in(1,\alpha)$
and $\delta':=\delta/(\delta-1)$,
 \beqlb\label{2.3}
 \ar\ar
\mbf{E}\Big\{\sup_{{\bf a}\in{\bf \bar{A}}}|H(Y_{\va,n,k}({\bf
a}))B_{\va,n,k} -H(\bar{Y}_{n,k}({\bf a},t,\va))B_t|
1_{U_{\va,\varsigma,\zeta}}\Big\} \cr
 \ar\ar\quad\le
C\Big\{\mbf{E}\big[|K_{\va,n,k}-K_{n,k}|^\delta 1_{U_{\va,\varsigma,\zeta}}\big]\Big\}^{1/\delta} \cr
 \ar\ar\qquad
+C\Big\{\mbf{E}\Big[\big||\tilde{M}_{\va,n,k,1}(t)|
+|\tilde{M}_{\va,n,k,2}(t)|+
|\tilde{N}_{\va,n,k}(t)|\big|^\delta1_{U_{\va,\varsigma,\zeta}}\Big]\Big\}^{1/{\delta'}}
 \eeqlb
for $t\in(t_{k-1},t_k]$ and
 \beqlb\label{2.13}
 \ar\ar
\mbf{E}\Big\{|H(Y_{\va,n,k}({\bf
\bar{a}}))B_{\va,n,k}
-H(K_{n,k})B_{\va,n,k}|1_{U_{\va,\varsigma,\zeta}}\Big\} \cr
 \ar\ar\qquad\le
C\Big\{\mbf{E}\Big[\Big||\bar{M}_{\va,n,k}(\bar{
a}_2)|^\delta+|K_{\va,n,k}-K_{n,k}|^\delta\Big| 1_{U_{\va,\varsigma,\zeta}}\Big] \Big\}^{1/\delta},
 \eeqlb
where
\begin{gather}
\tilde{N}_{\va,n,k}(t)
 :=
N_{\va,n,k} -y_0(t)^{-\frac1q},~\tilde{M}_{\va,n,k,1}(t) =
n\int_{t_{k-1}}^{t_k}y_\va([ns]/n)^{1-1/q}ds-y_0(t)^{1-1/q},\label{2.39}
\\ \tilde{M}_{\va,n,k,2}(t)
 :=
n\int_{t_{k-1}}^{t_k}y_\va(s)y_\va([ns]/n)^{-1/q}ds-y_0(t)^{1-1/q}
\label{2.40}.
 \end{gather}
For $n\ge2$ define functions $H_n(x)=H(x)1_{\{x>-n\}}$ and
$G_n(x)=H(x)1_{\{x\le-n\}}$. Then
$|H(x)-H(y)|\le |H_n(x)-H_n(y)|+|G_n(x)-G_n(y)|$ for $n\ge4$ and $x,y\in\mbb{R}$.
Let $\tilde{H}_n,\tilde{G}_n\in C^1(\mbb{R})$ satisfy
$\tilde{H}_n(x)=H_n(x)$ for all $x\in(-\infty,-n-1)\cup(-n,\infty)$
and $\tilde{G}_n(x)=G_n(x)$ for all
$x\in(-\infty,-n)\cup(-n+1,\infty)$. For large enough $n_1$,
$n_2:=n_1+2$ and $n_3:=n_1+4$, we have
$|H(x)-H(y)|\le
\sum_{k=1}^3[|\tilde{H}_{n_k}(x)-\tilde{H}_{n_k}(y)|
+|\tilde{G}_{n_k}(x)-\tilde{G}_{n_k}(y)|]$.
Then by \eqref{2.30} and the mean value theorem, there is a constant
$\tilde{c}_2=\tilde{c}_2(n_1)>0$ so that
$|\tilde{H}_{n_k}(x)-\tilde{H}_{n_k}(y)| \le \tilde{c}_2|x-y|$ and
 \beqnn
|\tilde{G}_{n_k}(x)-\tilde{G}_{n_k}(y)|
\le
|x-y|\int_0^1|\tilde{G}_{n_k}'(x+h(y-x))|dh
\le
\tilde{c}_2\big[|x|^\gamma1_{\{x<0\}}
+|y|^\gamma1_{\{y<0\}}\big]|x-y|
 \eeqnn
for $x,y\in\mbb{R}$.
It thus follows that
 \beqlb\label{2.23}
|H(x)-H(y)| \le \tilde{c}_2\big[1+|x|^\gamma1_{\{x<0\}}
+|y|^\gamma1_{\{y<0\}}\big]|x-y|, \qquad x,y\in\mbb{R}.
 \eeqlb
Since $1/\delta+1/\delta'=1$, by the H\"older inequality we have
 \beqlb\label{2.24}
 \ar\ar
\mbf{E}\Big\{\sup_{ {\bf a}\in{\bf \bar{A}}}|H(Y_{\va,n,k}({\bf a}))
-H(\bar{Y}_{n,k}({\bf a},t,\va))|1_{U_{\va,\varsigma,\zeta}}\Big\}\cr
 \ar\ar\le
\Big\{\mbf{E}\Big[\sup_{ {\bf a}\in{\bf \bar{A}}}
\big[1+|Y_{\va,n,k}({\bf a})|^\gamma 1_{\{Y_{\va,n,k}({\bf a})<0\}}
+|\bar{Y}_{n,k}({\bf a},t,\va)|^\gamma1_{\{\bar{Y}_{n,k}({\bf
a},t,\va)<0\}}\big]^{\delta'} 1_{U_{\va,\varsigma,\zeta}}\Big]\Big\}^{\frac{1}{\delta'}} \cr
 \ar\ar\quad
\cdot\Big\{\mbf{E}\big[\sup_{ {\bf a}\in{\bf \bar{A}}}
|Y_{\va,n,k}({\bf a})-\bar{Y}_{n,k}({\bf a},t,\va)|^\delta
1_{U_{\va,\varsigma,\zeta}}\big]\Big\}^{1/\delta},
\qquad t_{k-1}<t\le t_k.
 \eeqlb
Observe that
 \beqnn
\sup_{n,k\ge1,\,\va>0,\,t\in[0,1]}\big[|\tilde{M}_{\va,n,k,1}(t)|
+|\tilde{M}_{\va,n,k,2}(t)|
+|\tilde{N}_{\va,n,k}(t)|+|M_{\va,n,k}|+|N_{\va,n,k}|\big]
1_{U_{\va,\varsigma,\zeta}}<\infty.
 \eeqnn
Then by the fact $\delta'>\delta$ and H\"older inequality again we
get
 \beqlb\label{2.38}
 \ar\ar
\mbf{E}\Big\{|\bar{Y}_{n,k}({\bf a},t,\va)|\big[
|\tilde{N}_{\va,n,k}(t)| +|\tilde{M}_{\va,n,k,1}(t)|
+|\tilde{M}_{\va,n,k,2}(t)|\big]1_{U_{\va,\varsigma,\zeta}}\Big\} \cr
 \ar\le\ar
\Big\{\mbf{E}\big[|\bar{Y}_{n,k}({\bf
a},t,\va)|^\delta\big]\Big\}^{\frac{1}{\delta}}
\Big\{\mbf{E}\Big[\big[|\tilde{N}_{\va,n,k}(t)|
+|\tilde{M}_{\va,n,k,1}(t)|
+|\tilde{M}_{\va,n,k,2}(t)|\big]^{\delta'}1_{U_{\va,\varsigma,\zeta}}\Big]\Big\}^{\frac{1}{\delta'}} \cr
 \ar\le\ar
C\Big\{\mbf{E}\Big[\big[|\tilde{N}_{\va,n,k}(t)|
+|\tilde{M}_{\va,n,k,1}(t)|
+|\tilde{M}_{\va,n,k,2}(t)|\big]^{\delta}1_{U_{\va,\varsigma,\zeta}}\Big]\Big\}^{\frac{1}{\delta'}}.
 \eeqlb
Since $B\in C(\mbb{R}_+^2)$, $\sup_{n,k\ge1\,,\va>0}B_{\va,n,k}
1_{U_{\va,\varsigma,\zeta}}<\infty$.
In view of \eqref{2.23},
 \beqlb\label{2.36}
|H(x)|\le|H(x)-H(0)|+|H(0)|\le C\big(|x|^{\gamma+1}1_{\{x<0\}}+x1_{\{x>0\}}\big)+C,
 \eeqlb
which derives
 \beqnn
 \ar\ar
\big|H(Y_{\va,n,k}({\bf
a}))B_{\va,n,k}-H(\bar{Y}_{n,k}( {\bf
a},t,\va))B_t\big|\cr
 \ar=\ar
\big|H(\bar{Y}_{n,k}({\bf
a},t,\va))[B_{\va,n,k}-B_t]
+ B_{\va,n,k}\big[H(Y_{\va,n,k}({\bf
a}))-H(\bar{Y}_{n,k}({\bf a},t,\va))\big]\big|\cr
 \ar\le\ar
C\Big[1+|\bar{Y}_{n,k}({\bf
a},t,\va)|^{\gamma+1}1_{\{\bar{Y}_{n,k}({\bf
a},t,\va))<0\}}+\bar{Y}_{n,k}({\bf
a},t,\va)1_{\{\bar{Y}_{n,k}({\bf a},t,\va))>0\}}\Big]\cdot\Big[
|\tilde{N}_{\va,n,k}(t)| \cr
 \ar\ar
+[|a_2|+|\bar{a}_2|][|\tilde{M}_{\va,n,k,1}(t)|
+|\tilde{M}_{\va,n,k,2}(t)|]\Big]+C|H(Y_{\va,n,k}({\bf
a}))-H(\bar{Y}_{n,k}({\bf a},t,\va))|
 \eeqnn
on $U_{\va,\varsigma,\zeta}$ and
 \beqnn
|Y_{\va,n,k}({\bf a})-\bar{Y}_{n,k}({\bf a},t,\va)|
 \ar\le\ar
| a_3^{-1}[| a_2|+|\bar{
a}_2|]m_{\va,n}\big[|\tilde{M}_{\va,n,k,1}(t)|
+|\tilde{M}_{\va,n,k,2}(t)|\big] \cr
 \ar\ar
+|(\bar{ a}_1- a_1) a_3^{-1}|m_{\va,n}|\tilde{N}_{\va,n,k}(t)|
 + |\bar{ a}_3 a_3^{-1}||K_{\va,n,k}-K_{n,k}|.
 \eeqnn
Together with \eqref{2.2}--\eqref{2.1}, \eqref{2.24}--\eqref{2.38}
and Condition \ref{c1.2}(ii) one can get \eqref{2.3}.
By using \eqref{2.23},
 \beqnn
 \ar\ar
|H(Y_{\va,n,k}({\bf \bar{a}}))-H(K_{n,k})| \cr
 \ar\le\ar
C\Big[1+|Y_{\va,n,k}({\bf \bar{a}}))|^\gamma1_{\{Y_{\va,n,k}({\bf
\bar{a}}))<0\}}
+|K_{n,k}|^\gamma1_{\{K_{n,k}<0\}}\Big]\big[|\bar{M}_{\va,n,k}(\bar{
a}_2)|+|K_{\va,n,k}-K_{n,k}|\big].
 \eeqnn
Thus \eqref{2.13} follows from \eqref{2.2}--\eqref{2.1} and the
H\"older inequality.

{\bf Step 3.} Let $\delta''\ge1$, $\delta\in(1,\alpha)$ and $0\le
t'\le t''\le1$ with $|t'-t''|\le1/n$. Now
we show
 \beqlb\label{2.27}
\mbf{E}\Big\{|y_\va(t')^{1/\delta''}-y_0(t'')^{1/\delta''}|^\delta1_{U_{\va,\varsigma,\zeta}}\Big\} \le C[\va n^{-1/\alpha}+n^{-1}+\va]^\delta.
 \eeqlb
By using It\^o's formula on \eqref{1.1} one derives that
 \beqlb\label{2.12}
e^{\bar{ a}_2t}y_\va(t)=x_0+\bar{ a}_1\int_0^te^{\bar{ a}_2s}ds
+\va\bar{ a}_3\int_0^te^{\bar{
a}_2s}y_\va(s-)^{1/q}dz_0(s),\qquad t\in[0,1].
 \eeqlb
For $s\in[0,1]$ define $E_s=\{y_\va([ns]/n)\leq\varsigma,
y_\va(s-)\leq\varsigma\}$. Then on $U_{\va,\varsigma,\zeta}$,
 \beqnn
|y_\va(t')-y_\va(t'')|
 \ar\le\ar
e^{-\bar{ a}_2t''}|e^{\bar{ a}_2t'}y_\va(t')-e^{\bar{
a}_2t''}y_\va(t'')| +|1-e^{\bar{ a}_2(t'-t'')}|\varsigma \cr
 \ar\le\ar
|\va\bar{ a}_3|\Big|\int_{t'}^{t''}e^{\bar{
a}_2(s-t'')}y_\va(s-)^{1/q}1_{E_s}dz_0(s)\Big| +e^{|\bar{
a}_2|}(|\bar{a}_1|+\varsigma) n^{-1}.
 \eeqnn
It follows from \cite[Lemma 4.4]{LongQian} that for each $0<\delta_0<\alpha
$,
 \beqlb\label{2.10}
\mbf{E}\big\{|y_\va(t')-y_\va(t'')|^{\delta_0}1_{U_{\va,\varsigma,\zeta}}\big\}
 \le C\mbf{E}\Big\{\Big[ \int_{t'}^{t''}\va
y_\va(s)^{\frac{\alpha}{q}}1_{E_s}ds\Big]^{\frac{\delta_0}{\alpha}}\Big\} +\frac{C}{n^{\delta_0}}
\le C[\va
n^{-1/\alpha}+n^{-1}]^{\delta_0}.
 \eeqlb
By \eqref{2.12} and \cite[Lemma 4.4]{LongQian} again, for each
$t\in[0,1]$,
 \beqlb\label{2.25}
\mbf{E}\Big\{|y_\va(t)-y_0(t)|^\delta1_{U_{\va,\varsigma,\zeta}}\Big\} \le
\mbf{E}\Big[\Big|\va\bar{ a}_3\int_0^te^{\bar{
a}_2(s-t)}y_\va(s-)^{1/q}1_{E_s}dz_0(s)\Big|^\delta\Big]
 \le
C\va^\delta.
 \eeqlb
Observe that
$|z_1^{1/\delta''}-z_2^{1/\delta''}|\le C|z_1-z_2|$
for $z_1,z_2\in[\zeta,\varsigma]$.
Together with \eqref{2.10}--\eqref{2.25}
we have
 \beqnn
 \ar\ar
\mbf{E}\Big\{|y_\va(t')^{1/\delta''}-y_0(t'')^{1/\delta''}|^\delta1_{U_{\va,\varsigma,\zeta}}\Big\} \cr
 \ar\ar\quad\le
2^\delta\mbf{E}\Big\{|y_\va(t')^{1/\delta''}-y_\va(t'')^{1/\delta''}|^\delta1_{U_{\va,\varsigma,\zeta}}\Big\}
+2^\delta\mbf{E}\Big\{|y_\va(t'')^{1/\delta''}-y_0(t'')^{1/\delta''}|^\delta1_{U_{\va,\varsigma,\zeta}}\Big\} \cr
 \ar\ar\quad\le
C\mbf{E}\Big\{\Big[|y_\va(t')-y_\va(t'')|^\delta
+|y_\va(t'')-y_0(t'')|^\delta\Big]1_{U_{\va,\varsigma,\zeta}}\Big\}\le
C[\va n^{-1/\alpha}+n^{-1}+\va]^\delta,
 \eeqnn
which derives \eqref{2.27}.

{\bf Step 4.} Recall \eqref{1.3} and \eqref{2.39}--\eqref{2.40}. Let
$1<\delta<\alpha$. In this step we show for $t\in(t_{k-1},t_k]$,
\begin{gather}
\mbf{E}\Big[\Big||\tilde{M}_{\va,n,k,1}(t)|+|\tilde{M}_{\va,n,k,2}(t)|
+|\tilde{N}_{\va,n,k}(t)|\Big|^\delta1_{U_{\va,\varsigma,\zeta}}\Big]
 \le
C[\va n^{-1/\alpha}+n^{-1}+\va]^\delta,\label{2.11}\\
\mbf{E}\{|\bar{M}_{\va,n,k}(\bar{
a}_2)|^\delta1_{U_{\va,\varsigma,\zeta}}\} \le C[\va
n^{-1/\alpha}+n^{-1}]^\delta, \label{2.9} \\
\mbf{E}\{|K_{\va,n,k}-K_{n,k}|^\delta1_{U_{\va,\varsigma,\zeta}}\} \le C[\va n^{-1/\alpha}+n^{-1}]^{r\delta},\qquad r\in(0,1).
\label{2.4}
\end{gather}
Observe that on $U_{\va,\varsigma,\zeta}$, we have
$|\tilde{M}_{\va,n,k,1}(t)|^\delta
\le
n\int_{t_{k-1}}^{t_k}|y_\va([ns]/n)^{1-1/q}-y_0(t)^{1-1/q}|^\delta
ds$,
 \beqnn
|\tilde{M}_{\va,n,k,2}(t)|^\delta
 \ar\le\ar
\Big|n\int_{t_{k-1}}^{t_k}[y_\va(s)-y_\va([ns]/n)]y_\va([ns]/n)^{-1/q}ds
+\tilde{M}_{\va,n,k,1}(t)\Big|^\delta \cr
 \ar\le\ar
2^\delta
y_\va(t_{k-1})^{-\delta/q}n\int_{t_{k-1}}^{t_k}|y_\va(s)-y_\va([ns]/n)|^\delta
ds +2^\delta|\tilde{M}_{\va,n,k,1}(t)|^\delta \cr
 \ar\le\ar
2^\delta\zeta^{-\delta/q}n\int_{t_{k-1}}^{t_k}|y_\va(s)-y_\va([ns]/n)|^\delta
ds +2^\delta|\tilde{M}_{\va,n,k,1}(t)|^\delta, \cr
|\tilde{N}_{\va,n,k}(t)|^\delta
 \ar\le\ar
n\int_{t_{k-1}}^{t_k}\big|y_\va([ns]/n)^{-1/q}-y_0(t)^{-1/q}\big|^\delta
ds \cr
 \ar\le\ar
y_0(t)^{-\delta/q}\zeta^{-\delta/q}
n\int_{t_{k-1}}^{t_k}\big|y_\va([ns]/n)^{1/q}-y_0(t)^{1/q}\big|^\delta
ds
 \eeqnn
and
$|\bar{M}_{\va,n,k}(\bar{a}_2)|^\delta
\le
|\bar{ a}_2|^\delta\zeta^{-\delta/q}
n\int_{t_{k-1}}^{t_k}|y_\va([ns]/n)-y_\va(s)|^\delta ds$.
Then \eqref{2.11} and \eqref{2.9} follows from \eqref{2.27} and
\eqref{2.10}. Observe that for $t',t''\in[0,1]$,
 \beqnn
\big[y_\va(t')}{y_\va(t'')^{-1}\big]^{1/q}-1 =
y_\va(t'')^{-1-\frac1q}\Big\{\big[y_\va(t'')-y_\va(t')\big]y_\va(t')^{\frac1q}
-\big[y_\va(t'')^{1+\frac1q}-y_\va(t')^{1+\frac1q}\big]\Big\}.
 \eeqnn
Recall that $E_s=\{y_\va([ns]/n)\leq\varsigma,
y_\va(s-)\leq\varsigma\}$ for $s\in[0,1]$. Then on $U_{\va,\varsigma,\zeta}$,
 \beqnn
 \ar\ar
|K_{\va,n,k}-K_{n,k}|^\delta=\Big|n^{1/\alpha}\int_{t_{k-1}}^{t_k}
\Big[\big[y_\va(s-)}{y_\va([ns]/n)^{-1}\big]^{1/q}-1\Big]dz_0(s)\Big|^\delta
\cr
 \ar\ar\quad\le
2^\delta
y_\va(t_{k-1})^{-\delta-\frac{\delta}{q}}\Big|n^{1/\alpha}\int_{t_{k-1}}^{t_k}
\Big[\big[y_\va([ns]/n)-y_\va(s-)\big]y_\va(s-)^{\frac1q}\Big]dz_0(s)\Big|^\delta
 \cr
 \ar\ar\qquad+2^\delta
y_\va(t_{k-1})^{-\delta-\frac{\delta}{q}}\Big|n^{1/\alpha}\int_{t_{k-1}}^{t_k}
\big[y_\va([ns]/n)^{1+\frac1q}-y_\va(s-)^{1+\frac1q}\big]dz_0(s)\Big|^\delta
\cr
 \ar\ar\quad\le
2^\delta
\zeta^{-\delta-\frac{\delta}{q}}\Big|n^{1/\alpha}\int_{t_{k-1}}^{t_k}
\big[y_\va([ns]/n)-y_\va(s-)\big]y_\va(s-)^{\frac1q}1_{E_s}dz_0(s)\Big|^\delta
 \cr
 \ar\ar\qquad+2^\delta
\zeta^{-\delta-\frac{\delta}{q}}\Big|n^{1/\alpha}\int_{t_{k-1}}^{t_k}
\big[y_\va([ns]/n)^{1+\frac1q}-y_\va(s-)^{1+\frac1q}\big]1_{E_s}dz_0(s)\Big|^\delta.
 \eeqnn
Thus by \cite[Lemma 4.4]{LongQian} and Jensen's inequality again,
 \beqnn
 \ar\ar
\mbf{E}\Big\{|K_{\va,n,k}-K_{n,k}|^\delta1_{U_{\va,\varsigma,\zeta}}\Big\} \cr
 \ar\le\ar
C\mbf{E}\Big\{\Big[n\int_{t_{k-1}}^{t_k}\big[ |y_\va([ns]/n)
-y_\va(s)|^\alpha
y_\va(s)^{\frac{\alpha}{q}}+\big|y_\va([ns]/n)^{1+\frac1q}
-y_\va(s)^{1+\frac1q}\big|^\alpha\big]1_{E_s}ds\Big]^{\delta/\alpha}\Big\}.
 \eeqnn
Since
 \beqnn
|y_\va([ns]/n)^{1+\frac1q} -y_\va(s)^{1+\frac1q}|
\le
(1+1/q)\varsigma^{1/q}|y_\va([ns]/n) -y_\va(s)|
\le
C|y_\va([ns]/n) -y_\va(s)|^r
 \eeqnn
on $E_s$ by the mean value theorem, we know that
 \beqnn
\mbf{E}\{|K_{\va,n,k}-K_{n,k}|^\delta1_{U_{\va,\varsigma,\zeta}}\}
\le
C\{n\int_{t_{k-1}}^{t_k} \mbf{E}[|y_\va(s)
-y_\va([ns]/n)|^{r\alpha}1_{E_s}]ds\}^{\delta/\alpha}
 \eeqnn
by the H\"older inequality.
Combining with \eqref{2.10} we get \eqref{2.4}.

{\bf Step 5.} In this step we prove the following assertions:
\begin{gather}
n^{-1/2}\sum_{k=1}^n\big| [H(Y_{\va,n,k}({\bf \bar{a}}))
-H(K_{n,k})]B_{\va,n,k}\big|
\overset{p}{\longrightarrow}0,\label{2.14}\\
n^{-1}\sum_{k=1}^n \sup_{{\bf a}\in{\bf\bar{
A}}}\Big|H(Y_{\va,n,k}({\bf a}))B_{\va,n,k}
-\int_0^1H(\bar{Y}_{n,k}({\bf a},t,\va))B_tdt\Big|
\overset{p}{\longrightarrow}0.\label{2.6}
\end{gather}
It follows from \eqref{2.13} and \eqref{2.9}--\eqref{2.4} that for
any $r\in(\alpha/2,1)$
 \beqnn
\mbf{E}\Big\{\big| [H(Y_{\va,n,k}({\bf \bar{a}}))
-H(K_{n,k})]B_{\va,n,k}\big|1_{U_{\va,\varsigma,\zeta}}\Big\}\le C[\va n^{-1/\alpha}+n^{-1}]^r.
 \eeqnn
Observe that
 \beqnn
 \ar\ar
\Big| H(Y_{\va,n,k}({\bf a}))B_{\va,n,k}
-\int_0^1H(\bar{Y}_{n,k}({\bf a},t,\va))B_tdt\Big| \cr
 \ar\ar~\le
\sum_{i=1}^n\int_{t_{i-1}}^{t_i} |H(Y_{\va,n,k}({\bf
a}))B_{\va,n,k} -H(\bar{Y}_{n,k}({\bf
a},t,\va))B_t|dt.
 \eeqnn
Combining this with \eqref{2.3} and \eqref{2.11}--\eqref{2.4} we
obtain that for each $\delta\in(1,\alpha)$ and $r\in(0,1)$,
 \beqnn
\mbf{E}\Big\{\sup_{ {\bf a}\in{\bf \bar{ A}}}\Big|
H(Y_{\va,n,k}({\bf a}))B_{\va,n,k}
-\int_0^1H(\bar{Y}_{n,k}({\bf
a},t,\va))B_tdt\Big|1_{U_{\va,\varsigma,\zeta}}\Big\}\le
C[\va n^{-1/\alpha}+n^{-1}+\va]^{r(\delta-1)}.
 \eeqnn
It thus follows from
the Markov inequality that for each $\eta>0$ and $r\in(\alpha/2,1)$,
 \beqnn
 \ar\ar
\mbf{P}\Big\{n^{-1/2}\sum_{k=1}^n\big| [H(Y_{\va,n,k}({\bf
\bar{a}})) -H(K_{n,k})]B_{\va,n,k}\big|>\eta,
U_{\va,\varsigma,\zeta}\Big\} \cr
 \ar\ar\qquad\le
\eta^{-1}n^{-1/2}\sum_{k=1}^n\mbf{E}\Big\{\big| [H(Y_{\va,n,k}({\bf
\bar{a}}))
-H(K_{n,k})]B_{\va,n,k}\big|1_{U_{\va,\varsigma,\zeta}}\Big\} \cr
 \ar\ar\qquad\le
C\eta^{-1}[n^{1/2-r}+\va^r n^{(\alpha-2r)/(2\alpha)}] \to0
 \eeqnn
and
 \beqnn
 \ar\ar
\mbf{P}\Big\{\frac1n\sum_{k=1}^n \sup_{{\bf a}\in{\bf \bar{
A}}}\Big| H(Y_{\va,n,k}({\bf a}))B_{\va,n,k}
-\int_0^1H(\bar{Y}_{n,k}({\bf a},t,\va))B_tdt\Big|>\eta,
U_{\va,\varsigma,\zeta}\Big\} \cr
 \ar\ar\qquad\le
\eta^{-1}n^{-1}\sum_{k=1}^n\mbf{E}\Big\{\sup_{ {\bf a}\in{\bf \bar{
A}}}\Big| H(Y_{\va,n,k}({\bf a}))B_{\va,n,k}
-\int_0^1H(\bar{Y}_{n,k}({\bf
a},t,\va))B_tdt\Big|1_{U_{\va,\varsigma,\zeta}}\Big\}
\cr
 \ar\ar\qquad\le
C\eta^{-1}[\va n^{-1/\alpha}+n^{-1}+\va]^{r(\delta-1)}\to0
 \eeqnn
as $n\to\infty$. Then \eqref{2.14} and \eqref{2.6}
follow from \eqref{2.7}.

{\bf Step 6.} In this step we prove \eqref{2.34}. Since
$\bar{Y}_{n,k}( {\bf a},t,\va)-Y_{n,k}( {\bf
a},t)=(m_{\va,n}-m_0)Y({\bf a},t)$ and $H$ is uniformly continuous
on any bounded interval, we have as $n\to\infty$,
 \beqlb\label{2.28}
\sup_{{\bf a}\in{\bf \bar{ A}}}\int_0^1\big|[H(\bar{Y}_{n,k}({\bf
a},t,\va)) -H(Y_{n,k}({\bf a},t))]B_t\big|dt\to0.
 \eeqlb
By using \eqref{2.1} and
\eqref{2.36}, we get
$\mbf{E}\{\sup_{ {\bf a}\in{\bf \bar{ A}}}\int_0^1|H(Y_{n,k}(
{\bf a},t))B_t|dt\} <\infty$.
It is easy to see that for each $n,k\ge1$ and ${\bf a}\in{\bf \bar{ A}}$,
we have
$\mbf{E}\{\int_0^1H(Y_{n,k}({\bf a},t))B_tdt\} =\int_0^1
\<H(Y_0({\bf a},t,\cdot)),p\>B_tdt$.
Then by the uniform laws of large numbers (see e.g. \cite[Theorem
4.2.1]{Amemiya}),
 \beqnn
\sup_{{\bf a}\in{\bf \bar{ A}}}|n^{-1}\sum_{k=1}^n
\int_0^1H(Y_{n,k}({\bf a},t))B_tdt -\int_0^1
\<H(Y_0({\bf a},t,\cdot)),p\>B_tdt|
\overset{p}{\longrightarrow}0
 \eeqnn
as $n\to\infty$,
Hence \eqref{2.34} follows from \eqref{2.6} and \eqref{2.28}.

{\bf Step 7.} In this step we show \eqref{2.15}. Applying
\eqref{2.34} we obtain
 \beqnn
\frac1n\sum_{k=1}^n
\big[H(K_{n,k})B_{\va,n,k}+\bar{H}(K_{n,k})\big]^2
\overset{p}{\longrightarrow}\eta_0^2
 \eeqnn
as $n\to\infty$. Then by using \cite[Theorem
3.4]{HallH80} one can see that
 \beqlb\label{2.35}
n^{-1/2}\sum_{k=1}^n[H(K_{n,k})B_{\va,n,k}+\bar{H}(K_{n,k})]
\overset{d}{\longrightarrow} N(0,\eta_0^2)
 \eeqlb
as $n\to\infty$. Applying \eqref{2.14} we obtain
 \beqnn
n^{-1/2}\sum_{k=1}^n\Big|\big[H(Y_{\va,n,k}({\bf
\bar{a}}))B_{\va,n,k}+\bar{H}(Y_{\va,n,k}({\bf
\bar{a}}))\big]
-\big[H(K_{n,k})B_{\va,n,k}+\bar{H}(K_{n,k})\big]\Big|
\overset{p}{\longrightarrow}0
 \eeqnn
as $n\to\infty$. This, together with \eqref{2.35}, implies \eqref{2.15}, which completes the proof.
 \qed

{\bf Acknowledgement.}
The author is deeply grateful to Professors Zenghu Li and Chunhua Ma for their
encouragements and helpful discussions.
The author also thanks the referee and the editor for a number of suggestions on
the presentation of the results.

\end{document}